\documentclass[reqno,11pt]{amsart}
\usepackage{amssymb}
\usepackage{amsmath}
\usepackage{amsthm}
\usepackage{pgf}
\usepackage{color}
\usepackage{mathrsfs}
\usepackage{graphicx}   
\usepackage{hyperref}
\usepackage{graphicx}
\usepackage{tikz}
\usepackage{pgf}
\usepackage{enumerate}
\usepackage{bbm}
\usepackage{mathtools}
\usepackage{ulem}
\usetikzlibrary{positioning}
\usetikzlibrary{decorations.pathreplacing,decorations.markings}
\tikzset{->-/.style={decoration={
			markings,
			mark=at position #1 with {\arrow{>}}},postaction={decorate}}}
\usepackage{enumitem} 
\usepackage{comment}

\newtheorem{theorem}{Theorem}[section]
\newtheorem{lemma}[theorem]{Lemma}
\newtheorem{proposition}[theorem]{Proposition}
\newtheorem{corollary}[theorem]{Corollary}

\newcommand{\Rz}{\mathbb{R}}

\newcommand{\epsi}{\varepsilon}
\newcommand{\eps}{\varepsilon}

\def\eps{\varepsilon}

\renewcommand{\d}{{\mathrm d}}

\newcommand{\dt}{{\mathrm d} t}

\def\to{\rightarrow}

\newcommand{\weakto}{\rightharpoonup}

\numberwithin{equation}{section}

\makeatletter
\newcommand*\bigcdot{\mathpalette\bigcdot@{.5}}
\newcommand*\bigcdot@[2]{\mathbin{\vcenter{\hbox{\scalebox{#2}{$\m@th#1\bullet$}}}}}
\makeatother

\begin{document}
\title[WED approach for gradient flow in optimal control theory]
{Optimal control of gradient flows via\\ the Weighted Energy-Dissipation method}

\author{Takeshi Fukao}
  \address[Takeshi Fukao]{Faculty of Advanced Science and Technology,
  Ryukoku University, 1-5 Yokotani, Seta Oe-cho, Otsu-shi, Shiga 520-2194, Japan.}
\email{fukao@math.ryukoku.ac.jp}
\urladdr{https://fukao.math.ryukoku.ac.jp/indexe.html}

\author{Ulisse Stefanelli}
\address[Ulisse Stefanelli]{Faculty of Mathematics,  University of Vienna, Oskar-Morgenstern-Platz 1, A-1090 Vienna, Austria,
Vienna Research Platform on Accelerating Photoreaction Discovery,
University of Vienna, W\"ahringer Str. 17, 1090 Vienna, Austria, $\&$ Istituto di Matematica Applicata e Tecnologie Informatiche ``E. Magenes'' - CNR, v. Ferrata 1, I-27100 Pavia, Italy}
\email{ulisse.stefanelli@univie.ac.at}
\urladdr{http://www.mat.univie.ac.at/$\sim$stefanelli}

\author{Riccardo Voso}
\address[Riccardo Voso]{Faculty of Mathematics, University of Vienna, 
	Oskar-Morgenstern-Platz 1, A-1090 Vienna, Austria $\&$ Vienna School of Mathematics, Oskar-Morgenstern-Platz 1, A-1090 Vienna, Austria.}
\email{riccardo.voso@univie.ac.at}



\begin{abstract}
  We consider a general optimal control problem in the setting of
  gradient flows. Two approximations of the problem are presented, both
  relying on the variational reformulation of gradient-flow dynamics
  via the  Weighted-Energy-Dissipation variational approach. This
  consists in the minimization of global-in-time functionals over
  trajectories, combined with a limit passage. We show that the
  original nonpenalized problem and the two successive approximations
  admits solutions. Moreover, resorting to a $\Gamma$-convergence
  analysis we show that penalised optimal controls
  converge to nonpenalized one as the approximation is removed.
\end{abstract}

\maketitle

\section{Introduction}

This paper is concerned with an optimal control problem for abstract
gradient flows in Hilbert spaces.  We  are interested in
finding a solution to the following problem
\begin{equation}
\label{initialproblem}
\min_{(y,u)\in H^1(0,T;H)\times U}  P(y,u) \tag{$\mathcal{P}$}.
\end{equation}
The {\it control} $u:[0,T] \to H$ and the gradient flow
$y:[0,T]\to H$ are trajectories in the Hilbert space $H$ and $T>0$ is
some final reference time. The set $U$
of {\it admissible controls} is assumed to be compact in $L^2(0,T;H)$
and the functional $P$ is prescribed as by
$$ P(y,u) =
\begin{cases}
J(y,u)\quad &\text{if}\;y=S(u),\\
\infty &\text{else},
\end{cases} $$
where $J$ is a given {\it target functional}, which is assumed to  be lower semicontinuous with respect to the weak $\times$ strong
topology in $H^1(0,T;H)\times L^2(0,T;H)$. Moreover, $S(u)$ indicates the
unique solution $y\in H^1(0,T;H)$, given $u\in L^2(0,T;H)$ to the
{\it gradient flow problem}
\begin{align}
&\dot{y}(t) + \partial \phi(y(t))\ni  u(t) \quad  \text{for a.e.}\;t\in(0,T),\label{gradientflowi}\\
&y(0)=y^0.\label{gradientflowf}
\end{align}
The dot in \eqref{gradientflowi} indicates the time derivative of $y$. The functional $\phi : H
\rightarrow (-\infty,\infty]$ is assumed to be proper, lower
semicontinuous, and $\kappa$-convex for some $\kappa \in \mathbb{R}$,
i.e., $y \mapsto \psi(u)=\phi(y) - \kappa \lVert y \rVert^2/2$ is
convex. In particular, $\partial\phi(u) = \partial \psi(u)+\kappa u$, where
$\partial \psi(u)$ is the subdifferential in the sense of convex
analysis. Finally, $y^0\in D(\partial\phi)$.

The optimal control Problem \eqref{initialproblem} can be readily
proved to admit a solution, namely an optimal pair $(y,u)$. The focus of this note is on the possible
approximation of Problem \eqref{initialproblem}  by means of two
penalized problems. The departing point for such approximation is the
so called {\it Weighted Energy-Dissipation} (WED) approach to the
gradient flow problem. This consists in the
minimization of a family of $\epsi$-dependent 
global-in-time {\it WED} functionals $W_{\eps} : H^1(0,T;H) \times L^2(0,T;H)
\rightarrow (-\infty , \infty]$ given by 
\begin{equation}
\label{WEDfunctional}
W_{\eps}(y,u) =
\begin{cases}
\displaystyle\int_0^{T} e^{-t/\eps} \biggl( \frac{\eps}{2} \lVert \dot{y}(t)\rVert^2 + \phi(y(t)) - (u,y)  \biggr) \;dt  & \text{if}\; (y,u)\in K(y^0)\times U,\\
\qquad\qquad\qquad\qquad \infty  & \text{else},
\end{cases}
\end{equation}
where $\eps>0$ and the convex closed set $K(y^0)$ is given by $$K(y^0)
= \{ y \in H^1(0,T;H) : y(0) = y^0 \in H ,\, \phi \circ y \in L^1(0,T)
\}.$$
Indeed, given $u \in U$, the link between the minimization of
$W_\epsi(\cdot,u)$ and the gradient flow problem is revealed by
computing the corresponding Euler-Lagrange problem
\begin{align}
-\eps\ddot{y}_{\eps}(t) + \dot{y}_{\eps}(t) +\partial
  \phi(y_{\eps}(t)) &\ni u(t) \quad \text{for a.e.}\
  t\in(0,T), \label{eulerolagrangeepsilon}  \\
y_{\eps}(0) &= y^0\label{eulerolagrangeepsilon1},\\
\eps\dot{y}_{\eps}(T) &=0 \label{eulerolagrangeepsilon3}.
\end{align}
The latter is nothing but an elliptic-in-time regularization of the
gradient flow problem, which is recovered by taking $\epsi \to
0$. More precisely, owing to \cite{Mielke2}, for $\epsi$ small enough
one can uniquely define  
$$y^u_\epsi = \arg\min W_\epsi(\cdot,u)$$
and prove that $y^u_\epsi
\weakto y=S(u)$ in $H^1(0,T;H)$, see Proposition \ref{ms13}
below. Note the
occurrence of the final condition \eqref{eulerolagrangeepsilon3},
which eventually is dropped in the $\epsi \to 0$ limit.

In the setting of gradient flows, the WED approach has been applied to
mean curvature \cite{Ilmanen,spadarostefanelli}, to periodic solvability
\cite{Hirano94}, to micro-structure evolution \cite{ContiOrtiz}, to
the incompressible Navier-Stokes system
\cite{bathorystefanelli,ortizschmidtstefanelli}, and to stochastic PDEs
\cite{scarpastefanelli1}. The linear case is
mentioned in the
classical PDE textbook by Evans \cite[Problem~3, p.~487]{Evans}.

The general theory for $\kappa$-convex energies $\phi$ can
be found in 
\cite{Mielke2}, while the nonconvex setting is discussed in \cite{Akagi4}
and \cite{Melchionna201sei} deals with nonpotential perturbations. A
stability result via $\Gamma$-convergence is in
\cite{LieroMelchionna}  and
  \cite{melchionna2} uses the WED approach for studying  
symmetries of solutions. In the more general setting of metric spaces,
curves of maximal
slope 
have also been studied 
\cite{rossisavaresegattistefanelli1,rossisavaresegattistefanelli2}.

Besides the gradient flow case, the WED approach has also been
considered in the doubly nonlinear setting, see
\cite{MielkeOrtiz,Mielke1} for
rate-independent and
\cite{AkagiMelchionna,AkagiMelchionnaStefanelli,Akagi1,Akagi2, Akagi3}
for rate-dependent theories. In addition, the hyperbolic case of
semilinear waves has been considered in \cite{serratilli1,stefanelli},
also including  forcings
\cite{Mainini,tentarellitilli3,tentarellitilli,tentarellitilli2} or
dissipative terms  \cite{marveggio,Liero2,serratilli2}.
Applications to dynamic fracture
\cite{deLucaDalMaso,LarsenOrtizRichardson}, dynamic plasticity
\cite{DavoliStefanelli}, and Lagrangian mechanics \cite{Liero1} 
 are also available.

The aim of this paper is that of investigating two
approximations of Problem \eqref{initialproblem}, based on the WED functionals $W_{\eps}$. At first, we approximate the constraint $y = S(u)$ by considering the optimal control problem
\begin{equation}
	\label{epsilonproblem}
	\min_{(y,u)\in  H^1(0,T;H)\times U} P_{\eps}(y,u), \tag{$\mathcal{P}_{\eps}$} 
\end{equation}
where the approximating functional $P_\epsi$ is defined for all $\epsi
>0$ as
\begin{equation*}
	P_{\eps}(y,u) =
	\begin{cases}
	J(y,u)\quad &\text{if}\;y\in \arg\min W_{\eps}(\cdot,u),\\
	\infty &\text{else}.
	\end{cases}
\end{equation*}
This essentially amounts to replacing the gradient-flow constraint
\eqref{gradientflowi}-\eqref{gradientflowf} by its elliptic
regularization
\eqref{eulerolagrangeepsilon}-\eqref{eulerolagrangeepsilon3}. Problem
\eqref{epsilonproblem} is still a constrained
optimization. Nevertheless the constraint is itself expressed by a
minimization, turning \eqref{epsilonproblem} into a {\it bilevel}
optimization problem. The internal minimization layer is
actually a convex problem, which is therefore accessible to
efficient optimization techniques. Another  distintive feature of this
approach is that WED minimizers $y^u_\epsi$ turn out to be more
regular than the corresponding gradient flows $y \in
S(u)$, see Proposition \ref{ms13}.

As a second option, we consider a penalization of the constraint $y \in \arg\min
W_\epsi(\cdot,u)$ in \eqref{epsilonproblem}, tuning to an unconstrained optimal control problem, namely,
\begin{equation}
\label{epsilonlambdaproblem}
\min   P_{\eps\lambda}(y,u)  \tag{$\mathcal{P}_{\eps\lambda}$},
\end{equation}
with $P_{\eps\lambda}$ given for all $\epsi,\, \lambda >0$ as
\begin{equation*}
P_{\eps\lambda}(y,u) = 
J(y,u) + \frac{1}{\lambda}\bigl(W_{\eps}(y,u) - M_\epsi^u  \bigr).
\end{equation*}
Here, $\lambda>0$ serves as penalization parameter and $M_\epsi^u$ denotes the minimum value of $W_{\eps}(\cdot,u)$,
given $u$, i.e., $M_\epsi^u = W_{\eps}(y^u_\epsi,u)$. In particular,
the functional $y\mapsto W_{\eps}(y,u) - M_\epsi^u $ is nonnegative
and vanishes iff $y=y^u_\epsi$. Although the unconstrained formulation
of \eqref{epsilonlambdaproblem} is appealing, one shall stress that
the computation of $y^u_\epsi$, i.e., the minimization of
$W_\epsi(\cdot,u)$ is still needed for evaluating the minimum value
$M^u_\epsi$, see Lemma \ref{lemma1} below.
A remarkable feature,
however, is that the functional $y\mapsto W_{\eps}(y,u) - M_\epsi^u $
controls the distance $\| y - y^u_\epsi\|^2_{H^1(0,T;H)}$, which may
be of some applicative interest.

Our main result is Theorem \ref{theorem} below. We first check that Problems \eqref{initialproblem},
\eqref{epsilonproblem}, and \eqref{epsilonlambdaproblem} are solvable,
namely, there exist optimal pairs $(y,u)$, $(y_\epsi,u_\epsi)$, and
$(y_{\epsi\lambda},u_{\epsi\lambda})$, respectively (Theorem
\ref{theorem}.i). Then, we prove that $$ P_{\eps\lambda}
\stackrel{\Gamma}{\to} P_{\eps}\quad \text{and}\quad  P_{\eps}
\stackrel{\Gamma}{\to} P $$ in the sense of
$\Gamma$-convergence with respect to the weak$\times$strong topology
of  $H^1(0,T;H) \times L^2(0,T;H)$, which we indicate as $\tau$. Upon checking the coercivity of
$P_\epsi$ and $P_{\epsi\lambda}$ with respect to topology $\tau$,
this allows us to prove that, for some not relabeled subsequences,
$$(y_\epsi,u_\epsi)\stackrel{\tau}{\to} (y,u) \quad \text{and} \quad
(y_{\epsi\lambda},u_{\epsi\lambda}) \stackrel{\tau}{\to} (y_\epsi,u_\epsi)$$
as $\epsi \to 0$ and $\lambda\to 0$, respectively, (Theorem
\ref{theorem}.ii-iii).

Eventually, we tackle the joint $(\epsi,\lambda)\to (0,0)$ in the
specific case $\lambda = \lambda_\epsi$ with $\limsup_{\eps \rightarrow 0}
\lambda_\eps \eps^{-3} e^{T/\eps} = 0$ proving that
$$P_{\eps\lambda} \stackrel{\Gamma}{\rightarrow} P.$$
Again due to coerciveness, this entails that $(y_{\epsi,\lambda})
\stackrel{\tau}{\to} (y,u)$ along some not relabeled subsequence
(Theorem \ref{theorem}.iv).

The plan of the paper is the following. In Section \ref{notation}, we
introduce the assumptions and state our main result, namely, Theorem
\ref{theorem}. Section \ref{lemmas} collects a series of lemmas, which
will be used throughout. The proof of Theorem \ref{theorem} is in
Sections \ref{well}-\ref{proof}. Specifically, existence for the Problems
\eqref{initialproblem}, \eqref{epsilonproblem}, and
\eqref{epsilonlambdaproblem} is checked in Section
\ref{well} and convergences as $\epsi \to 0$ and
$\lambda\to0$ are discussed in Section \ref{proof}.

\section{Statement of main results}
\label{notation}
In this section, we introduce assumptions, recall some result from
\cite{Mielke2}, and state our main result, i.e., Theorem
\ref{theorem}. Let us start by assuming that
\begin{itemize}
	\item[\textbf{(A1)}] $H$ is a real Hilbert space,
          $T>0$, and
          $U \subset L^2(0,T;H)$ is nonempty and compact.
          \end{itemize}
We indicate by $\lVert \cdot \rVert$ and by $(\cdot,\cdot)$ the norm
and the inner product in $H$, respectively, and we indicate by $\lVert
\cdot \rVert_E$ the norm in the generic Banach space $E$.

Concerning the functional $\phi$ we ask that
\begin{itemize}
	\item[\textbf{(A2)}] $\phi : H
	\rightarrow (-\infty,\infty]$ is proper, $\kappa$-convex, and
        lower semicontinuous, and $y^0 \in D(\partial \phi)$.
      \end{itemize}
      In particular, the {\it effective domain} $D(\phi)=\{v\in H
      \::\:\phi(v) < \infty\}$ is not empty. Moreover,
$$ v \mapsto \psi(v) = \phi(v) - \frac{\kappa}{2} \lVert v
\rVert^2\quad \text{is convex}, $$ and $D(\psi)  =
D(\phi)$. Equivalently, one can state $\kappa$-convexity of $\phi$ as $$ \phi(rw +
(1-r)v) \leq r \phi(w) + (1 - r) \phi(v) - \frac{\kappa}{2}r(1-r)
\lVert w - v \rVert^2 \quad \forall w,v \in H,\; 0\leq r \leq 1. $$
Correspondingly, we define the (Fr\'echet) subdifferential $\partial
\phi : H \to 2^H$ as $\partial \phi (y) = \partial \psi(y)+\kappa y$,
where $\partial \psi $ is the subdifferential of $\psi$ in the sense
of convex analysis \cite{brezis}. This implies that one has $\eta \in \partial
\phi(y)$ iff $y \in D(\phi)$ and
$$(\eta,w-y) \leq \phi(w) - \phi(y) + \frac{\kappa}{2}\| w-y\|^2\quad
\forall w \in D(\phi),$$
$D(\partial \phi)=\{y\in D(\phi) \::\: \partial \phi(y)\not =
\emptyset\}=D(\partial \psi)$, and $\partial \phi(y)$ is convex and
closed, for all $y \in D(\partial \phi)$. In particular, for all  $y
\in D(\partial \phi)$ one has that $\partial
\phi(y)$ has a unique element of minimal norm, which we indicate by
$(\partial \phi(y))^\circ$.

The next proposition summarizes results from \cite{Mielke2}.
\begin{proposition}[WED approach to gradient flows]
	\label{ms13}
	Under assumptions {\rm (A1)-(A2)}, there exists $\eps_0 >0$ so
        that for all $\eps \in (0,\eps_0)$ and all $u \in U$ the
        functional $W_{\eps}(\cdot,u)$ is $\kappa_{\eps}$-convex in
        $H^1(0,T;H)$ for some $\kappa_\epsi>0$. In particular, there
        exists a unique minimizer $y_{\eps}^u = \arg\min W_{\eps}
        (\cdot,u)$. One has that $y^{u}_{\eps} \in H^2(0,T;H)$ is the
        unique solution of the Euler-Lagrange problem
        \eqref{eulerolagrangeepsilon}-\eqref{eulerolagrangeepsilon3}
        and that $\eta_{\eps}^{u} = u +\epsi \ddot{y}^u_\epsi
        -\dot{y}^u_\eps  \in L^2(0,T;H)$ and fulfills $\eta^u_\epsi
        \in \partial \phi(y^u_\epsi)$ a.e.
	Moreover, there exists a nondecreasing function $\ell : \mathbb{R}_+ \rightarrow \mathbb{R}_+$ independent of $\eps$ such that
	\begin{align}
	\label{maxreg}
		&\eps \lVert \ddot{y}^{u}_{\eps} \rVert_{L^2(0,T;H)} +
           \eps^{1/2} \lVert \dot{y}^{u}_{\eps}
           \rVert_{L^{\infty}(0,T;H)} + \lVert
           \dot{y}^{u}_{\eps}\rVert_{L^2(0,T;H)} + \lVert
           \eta_{\eps}^{u} \rVert_{L^2(0,T;H)}\nonumber\\ &\quad \leq \ell (\lVert u \rVert_{L^2(0,T;H)} + \lVert y^0 \rVert + \lVert (\partial \phi (y^0))^{\circ}\rVert),
	\end{align}
	where $(\partial \phi (y^0))^{\circ}$ is the element of
        minimal norm in $\partial \phi (y^0)$.	As $\epsi \to 0$, one
        has that  $y_{\eps}^{u} $ converges to $y \in S(u)$ weakly in
        $H^1(0,T;H)$ and strongly in $H^\sigma (0,T;H)$ for all
        $\sigma \in (0,1)$. 
      \end{proposition}
      
In all of the following, we will tacitly assume that
$\eps \in (0,\eps_0)$, so that Proposition \ref{ms13} holds. In
particular, $ y^u_\epsi = \arg\min W_\epsi(\cdot ,  u)$ is well-defined. 

Concerning the target functional $J$ we assume that 
\begin{itemize}
	\item[\textbf{(A3)}] $J:H^1(0,T;H)\times L^2(0,T;H) \rightarrow \mathbb{R}_+$, $\mathbb{R_+} := [0,\infty)$, is lower semicontinuous with respect to $\tau$.
	Moreover, for all $u \in L^2(0,T;H)$ given, $y \mapsto J(y,u)$
        is upper semicontinuous with respect to the strong
        $H^{\sigma}(0,T;H)$ topology, for some $\sigma \in (0,1)$.
      \end{itemize}
      A possible example of functional $J$ fulfilling (A3) is
      $$J(y,u) = f(y(T))+ \int_0^Tg(y,u)\, {\rm d } t  $$
      where $f:H\to
      \Rz$ is continuous and $g:H \times H \to \Rz$ is continuous and
      bounded.

      Let us recall that, given functionals $F_{\rho},\, F :H^1(0,T;H)
      \times L^2(0,T;H) \to \Rz\cup\{\infty\}$ for $\rho>0$, we say
      that the
      sequence $(F_\rho)_\rho$ 
      {\it $\Gamma$-converges} to  $F$ with respect to the topology
      $\tau$ and we write $F = \Gamma_\tau\lim_{\rho \to 0} F_\rho$  if the following conditions hold
\begin{enumerate}
	\item[{(\rm i)}] ($\Gamma$-$\liminf$ inequality) For every
          $(y_{\rho},u_{\rho}) \overset{\tau}{\rightarrow} (y,u)$ we have
	\begin{equation}
	\label{liminf1}
	F(y,u)\leq\liminf_{\rho\to 0} F_{\rho}(y_{\rho},u_{\rho}) ;
	\end{equation}
	\item[{(\rm ii)}] (Recovery sequence) For every $(\hat{y},\hat{u}) \in H^1(0,T;H)\times L^2(0,T;H)$ there exists $(\hat{y}_{\rho},\hat{u}_{\rho}) \overset{\tau}{\rightarrow} (\hat{y},\hat{u})$ such that
	\begin{equation}
	\label{limsup1}
	\limsup_{ \rho\to 0} F_{\rho}(\hat{y}_{\rho},\hat{u}_{\rho}) \leq F(\hat{y},\hat{u}).
	\end{equation}
\end{enumerate}

      We now state our main result.
\begin{theorem}[WED approach to optimal control]
	\label{theorem}
	Assume  {\rm (A1)-(A3)}. Then:
	\begin{enumerate}
		\item[{\rm i)}] For all $\eps,\lambda>0 $ Problems \eqref{initialproblem}, \eqref{epsilonproblem}, and \eqref{epsilonlambdaproblem} admit solutions.
		\item[{\rm ii)}] $P_\epsi \stackrel{\Gamma}{\to} P$ as
                  $\epsi \to 0$. Any sequence $(y_{\eps},u_{\eps})_{\eps}$ of solutions to Problem \eqref{epsilonproblem} admits a not relabeled subsequence such that $(y_{\eps},u_{\eps})\overset{\tau}{\rightarrow}(y,u)$ where $(y,u)$ solves Problem~\eqref{initialproblem}.
		\item[{\rm iii)}] $P_{\epsi\lambda} \stackrel{\Gamma}{\to} P_\epsi$ as
                  $\lambda \to 0$, for all $\eps>0$ fixed. Any sequence $(y_{\eps\lambda},u_{\eps\lambda})_\lambda$ of solutions to the Problem \eqref{epsilonlambdaproblem} admits a not relabeled subsequence such that $(y_{\eps\lambda},u_{\eps\lambda}) \overset{\tau}{\rightarrow} (y_{\eps},u_{\eps})$ where $(y_{\eps},u_{\eps})$ solves Problem~\eqref{epsilonproblem}.
		\item[{\rm iv)}] Let $\lambda =
                  \lambda_\epsi$ with $\limsup_{\epsi \to 0} \lambda_\epsi\eps^{-3} e^{-T/\eps} = 0$. Then
                  $P_{\epsi\lambda_\eps}\stackrel{\Gamma}{\to}P$ as $\epsi
                  \to 0$. Any sequence of solutions
                  $(y_{\epsi\lambda_\epsi},u_{\epsi\lambda_\eps})$ to  Problem
                 \eqref{epsilonlambdaproblem} admits a
                  not relabeled subsequence such that $(y_{\eps\lambda_\eps},u_{\eps\lambda_\eps}) \overset{\tau}{\rightarrow} (y,u)$ where $(y,u)$ solves Problem \eqref{initialproblem}.
	\end{enumerate}	
      \end{theorem}

  The proof of Theorem \ref{theorem} is given in Section
  \ref{well}-\ref{proof}. More precisely, we give a proof of Theorem \ref{theorem}.i in Section \ref{well} whereas Theorem \ref{theorem}.ii-iv is
  proved in Section \ref{proof}.

  \subsection{An example}
  Before closing this section, we give an illustration of the results
  by resorting to simple ODE example. In particular, we consider the
  ODE 
  \begin{equation}\label{eq:ODE} \dot y + y = u, \quad y(0)=1
    \end{equation}
  with $U=\{u \::\: u(t)= u_0 e^{-t} \ \text{for some} \ u_0 \in
  [0,1]\}$. We are interested in minimizing
  $$J(y,u) = \frac{1}{2} \int_0^1 (y(t) - e^{-t})^2 \,\dt +
  \frac{1}{2} \int_0^1 t^2 (u(t) - e^{-t})^2 \, \dt.$$
  Note that this fits in the theory by letting $H=\Rz$, $T=1$, and
  $\phi(y)=y^2/2$.
  In particular, $U$ is clearly compact into $L^2(0,1)$.

  Problem $\mathcal P$ reads
  \begin{align*}
	\min_{u_0 \in [0,1]} \biggl\{ J(y,u) \ : \ \dot y(t) + y(t) =  u_0 e^{-t},   \ y(0) = 1 \biggr\},
  \end{align*}
  and can be directly solved. For all $u_0$ the solution of
  \eqref{eq:ODE} is $y(t)= (1+tu_0) e^{-t}$ and we have
  \begin{align*}
  J(y,u) &=\frac12 \int_0^1  u_0^2  t^2
  e^{-2t}\, \d t  + \frac12 \int_0^1 (u_0-1)^2 t^2
  e^{-2t}\, \d t  \\
         & =  
           \frac{1}{8}(1-5e^{-2})  \big( u_0^2 + (u_0-1)^2\big)
           \end{align*}
which is minimized at $u_0=1/2$ with value $ (1 - 5 e^{-2})/16$ and a corresponding optimal
  trajectory $y(t) = (1 + t/2)  e^{-t}$.

  Let us now turn to Problem $\mathcal P_\epsi$ which can be written as
 \begin{align*}
	\min_{u_0 \in [0,1]} \biggl\{ J(y,u) \ : \ 
     -\epsi \ddot y(t)+\dot y(t)  + y(t) = u_0 e^{-t},   \ y(0) = 1, \ 
    y'(1)=0  \biggr\}.
 \end{align*}
 Given $u \in U$, the only solution $y^u_\epsi$ to
 $$-\epsi \ddot y(t)+\dot y(t)+ y(t) = u_0e^{-t}, \quad  y(0) = 1, \quad  y'(1) =
 0 $$
 is given by
\begin{equation*}
	y_{\eps}^u (t) = c_{\eps}^- e^{r_{\eps}^-t} + c_{\eps}^+ e^{r_{\eps}^+t} - \frac{u_0}{\eps} e^{-t},
\end{equation*}
where
\begin{align*}
	&r_{\eps}^- = \frac{ 1 - \sqrt{4 \eps + 1}}{2 \eps},\qquad r_{\eps}^+ = \frac{ 1 + \sqrt{4 \eps +1}}{2 \eps},\\
	&c_{\eps}^- = \frac{ \displaystyle \frac{u_0}{  \eps e} + r_{\eps}^+\bigl( 1 + \frac{u_0}{\eps} \bigr) e^{r_\eps^+}}{r_{\eps}^+ e^{r_{\eps}^+} - r_{\eps}^- e^{r_{\eps}^-}},\qquad c_{\eps}^+ = 1 + \frac{u_0}{\eps} - c_{\eps}^-.
\end{align*}
The value of $J$ at $(y^u_\eps,u)$ can be explicitly computed as
a function of $u_0$ as
\begin{align*}
  u_0 &\mapsto j_\epsi(u_0) :=	J(y_{\eps}^u,u_0 e^{-t}) \\
                          &=  \frac{1}{2} \Bigg(
                                    \frac{(c_{\eps}^-)^2 }{2
                                    r_{\eps}^-} \bigl(	e^{2
                                    r_{\eps}^-} -1	\bigr) +
                                    \frac{(c_{\eps}^+)^2}{2
                                    r_{\eps}^+} \bigl( e^{2r_{\eps}^+}
                                    -1 \bigr) 	- \frac{1}{2} \biggl(
                                    \frac{u_0}{\eps} + 1\biggr)^2
                                    (e^{-2} -1) \\&\quad+ \frac{2
  c_{\eps}^- c_{\eps}^+}{r_{\eps}^- + r_{\eps}^+} \bigl( e^{r_{\eps}^-
  + r_{\eps}^+} -1  \bigr)	- \frac{2 c_{\eps}^-}{r_{\eps}^--1}
  \biggl( \frac{u_0}{\eps} + 1\biggr) \bigl( e^{r_{\eps}^- -1} -1
  \bigr) \\&\quad - \frac{2 c_{\eps}^+}{r_{\eps}^+-1} \biggl(
  \frac{u_0}{\eps} +1 \biggr) \bigl( e^{r_{\eps}^+ -1} -1 \bigr)
  \Bigg)				 +\frac{1}{8} \bigl( 1 - 5
  e^{-2}	\bigr)  (u_0 -1)^2 
\end{align*} 
A tedious but elementary computation ensures that $u_0 \mapsto
J(y_{\eps}^u,u_0 e^{-t})$ converges to $u_0 \mapsto
(1/8)(1-5e^{-2})  \big( u_0^2 + (u_0-1)^2\big)$ uniformly on
$ [0,1]$.  In particular, the minimum value of $j_\eps$ on
$[0,1]$ converges to that of $j_0$, namely $ (1 - 5 e^{-2})/16$,  and
the minimizers of $j_\eps$ converge to the minimizer
$1/2$, as well.

We conclude by considering Problem $\mathcal P_{\eps \lambda}$. This
amount to the following
\begin{align*}
	\min_{\substack{u_0 \in [0,1], \\ y \in H^1(0, 1), \\ y(0)=1}} \biggl\{ J(y,u)
  +\frac{1}{\lambda} \left(\int_0^1e^{-t/\epsi}\left(\frac{\epsi}{2} |\dot
  y(t)|^2 +\frac12|y(t)|^2 - u_0 e^{-t}y(t)\right)\, \d t - M^u_\epsi
  \right)  \biggr\}.
\end{align*}
In order to tackle this minimization problem, one needs to evaluate
$M^u_\epsi$, which generally calls for another minimization. In this
example however one can use the above expression for  $y^u_\epsi$ and explicitly
compute $M^u_\epsi$ 
\begin{align*}
M_{\eps}^u =& \frac{\eps^2 (r_{\eps}^-)^2 (c_{\eps}^-)^2}{4\eps r_{\eps}^- - 2} \bigl( e^{2r_{\eps}^- - 1/{\eps}} -1	\bigr) + \frac{\eps^2 (r_{\eps}^+)^2 (c_{\eps}^+)^2}{4 \eps r_{\eps}^+ -2} \bigl( e^{2r_{\eps}^+ - 1/{\eps}} -1  \bigr) \frac{u_0^2}{ 4\eps + 2} \bigl( e^{-2 - 1/{\eps}} -1   \bigr) \\&+ \frac{c_{\eps}^- c_{\eps}^+ r_{\eps}^- r_{\eps}^+}{ \eps (r_{\eps}^- + r_{\eps}^+ ) -1} \bigl( e^{r_{\eps}^- + r_{\eps}^+ - 1/{\eps}} -1 \bigr) + \frac{\eps c_{\eps}^-}{\eps r_{\eps}^- -1} \bigl( e^{r_{\eps}^- - 1/{\eps}} -1  \bigr) \\&+ \frac{\eps c_{\eps}^+}{\eps r_{\eps}^+ -1} \bigl( e^{r_{\eps}^+ - 1/{\eps}} -1  \bigr) + \frac{u_0}{1+\eps } \bigl( e^{-1 - 1/{\eps}} -1  \bigr).
\end{align*}

      \section{Preliminary lemmas}
      \label{lemmas}

      Before moving to the proof of Theorem \ref{theorem}, we present
      in this section some preliminary lemmas, which will be used
      throughout and which complement the
      analysis in~\cite{Mielke2}.

    To start with, let us explicitly remark that, although $\phi$
      is just $\kappa$-convex, the classical tools from convex
      analysis apply to $\partial \phi$, as well. In particular, we
      have that
      \begin{align}
        \label{eq:limsup}
        y_n \weakto y, \ \eta_n \weakto \eta, \ \eta_n \in \partial
        \phi (y_n), \ \limsup_n (\eta_n,y_n) \leq (\eta,y) \
        \Rightarrow \ \eta \in \partial\phi(y).
      \end{align}
      Indeed, one has that $\eta_n -\kappa y_n \in \partial
      \psi(y_n)$, $\eta_n -\kappa y_n \weakto \eta-\kappa y$, and
      \begin{align*}
        \limsup_n (\eta_n - \kappa y_n,y_n) \leq \limsup_n
        (\eta_n,y_n) - \liminf_n \kappa \| y_n\|^2 \leq (\eta - \kappa
        y,y)
      \end{align*}
so that using  \cite[Prop~2.5, p.~27]{Brezis73}  one finds $\eta - \kappa y \in
\partial \psi(y)$, which entails $\eta\in \partial \phi(y) $. The
identification \eqref{eq:limsup} equivalently holds in its integrated form for sequences
$y_n$ and $\eta_n$ weakly converging in $L^2(0,T;H)$, namely.
\begin{align}
        &y_n \weakto y, \ \eta_n \weakto \eta \ \text{in} \ L^2
          (0,T;H),  \ \eta_n \in \partial
        \phi (y_n) \ \text{a.e.}, \nonumber\\
         &\qquad
            \limsup_n \int_0^T(\eta_n,y_n) \,
    \d t\leq\int_0^T (\eta,y) \, \d t\\
  &\quad 
        \Rightarrow \ \eta \in \partial\phi(y) \ \text{a.e.}
        \label{eq:limsup2}
      \end{align}

One can also readily prove the following generalizaton
        to the $\kappa$-convex case of the classical chain rule
        \cite[Lemme 3.3, p.~73]{brezis}
        \begin{align}
          \label{eq:chain}
   & y \in H^1(0,T;H), \ \eta \in L^2(0,T;H), \ \eta \in \partial
          \phi(u) \ \text{a.e. in} \ (0,T) \\
          &\quad\Rightarrow \ \phi \circ
            y \ \text{is absolutely continuous on $[0,T]$ and} \\
          &\qquad 
          \frac{\d }{\d t} \phi\circ y = (\eta,y) \ \text{a.e. in} \ (0,T).     
        \end{align}

In the following, we use the symbol $c$ to indicate a generic positive
constant, possibly depending on $T,\,U,\,y^0$ but independent on
$\epsi$ and $\lambda$ and possibly varying from line to
line.

      We are now ready to present the lemmas.

\begin{lemma}[Value of $M_\epsi^u$]\label{lemma1} For all $u\in U$,
  recalling that $y_{\eps}^u =
  \arg \min W_{\eps}(\cdot, u)$ and $M_\epsi^u = W_{\eps}(y_{\eps}^u ,
  u)$ we have
	\begin{align*}
	M_\epsi^u &= -\frac{\eps^2}{2} \lVert \dot{y}_{\eps}^u(0)
   \rVert^2 - \eps e^{-T/\eps} \phi(y^u_{\eps}(T)) + \eps \phi(y^0)
          \\
          &\quad +
   \eps \int_0^T e^{-t/\eps}(u,\dot{y}_{\eps}^u) \, \d t - \int_0^T e^{-t/\eps}(u,y_{\eps}^u) \,\d t.
	\end{align*}
\end{lemma}
\begin{proof}
	The trajectory $y_{\eps}^u$ solves the Euler-Lagrange problem
        \eqref{eulerolagrangeepsilon}-\eqref{eulerolagrangeepsilon3}. By
        taking the scalar product with $e^{-t/\eps}\dot{y}_{\eps}^u\in
        H^1(0,T;H)$  in equation \eqref{eulerolagrangeepsilon},
        integrating in time, and using the chain rule \eqref{eq:chain}, we get
	\begin{align*}
	\int_0^T e^{-t/\eps} \biggl(	-\frac{\eps}{2} \frac{\d}{\d t} \lVert \dot{y}_{\eps}^u \rVert^2 + \lVert \dot{y}_{\eps}^u \rVert^2 + \frac{\d}{\d t} \phi(y_{\eps}^u) - (u , \dot{y}_{\eps}^u)		\biggr) \, \d t = 0.
	\end{align*}
	By integrating by parts the first and the third term, one obtains
	\begin{align*}
	&\biggl[e^{-t/\eps}\biggl(
   -\frac{\eps}{2}\lVert\dot{y}_{\eps}^u\rVert^2 + \phi(y_{\eps}^u)
   \biggr)\biggr]_{0}^T -  \int_0^T e^{-t/\eps}(u,\dot{y}_{\eps}^u)\,
   \d t \\
          &\quad + \frac{1}{\eps} \int_0^T e^{-t/\eps} (u , y_{\eps}^u) \, \d t + \frac{1}{\eps}M_\epsi^u = 0,
	\end{align*}
	where we used the definition \eqref{WEDfunctional} of the WED functional.
	The thesis follows from conditions \eqref{eulerolagrangeepsilon1}-\eqref{eulerolagrangeepsilon3}. 
\end{proof}

\begin{lemma}[Continuity of the map $u \mapsto
  y_{\eps}^u$]\label{lemma2} Let $( u_k)_k \subset U$ be such that
  $u_k \rightarrow u$ in $L^2(0,T;H)$ and let $\eta_\eps^{u_k} = \eps \ddot{y}^{u_k}_\epsi - \dot{y}^{u_k} _\epsi+u_k$. Up to not relabeled subsequences, one has
	\begin{align*}
	&y_{\eps}^{u_k} \rightharpoonup y_{\eps}^u \quad \text{in} \quad H^2(0,T;H),\\
	&\eta_{\eps}^{u_k} \rightharpoonup \eta_{\eps}^u \quad \text{in}\quad L^2(0,T;H)
	\end{align*}
	where $\eta_\eps^u = \eps \ddot{y}^u_\epsi - \dot{y}^u_\epsi
        +u $.
\end{lemma}
\begin{proof}
	From the uniform estimate \eqref{maxreg} we may extract not
        relabeled subsequences such that $y_{\eps}^{u_k} \rightharpoonup y$ in $H^2(0,T;H)$ and $\eta_{\eps}^{u_k} \rightharpoonup \eta$ in $L^2(0,T;H)$ and get
	\begin{equation*}
	-\eps \ddot{y}(t) + \dot{y}(t) + \eta(t) = u(t) \quad \text{in}\;H\;\text{a.e.}\;t\in(0,T).
      \end{equation*}
      As $y^{u_k}_\epsi(t)\weakto y(t)$ and
      $\dot{y}^{u_k}_\epsi(t)\weakto \dot{y}(t)$ for all $t \in [0,T]$
      we have that $y(0)=y^0$ and $\eps \dot{y}(T)=0$. In order to
      conclude the proof it hence suffices to check that $\eta \in
      \partial \phi(y)$ a.e. Take the scalar 
      $\limsup$ of the integral over $(0,T)$ of the scalar product between $\eta_{\eps}^{u_k}$ and $y_{\eps}^{u_k}$. Using equation \eqref{eulerolagrangeepsilon} at level $k$, we obtain
	\begin{align*}
	&\limsup_{k\rightarrow \infty} \int_0^T
   (\eta_{\eps}^{u_k},y_{\eps}^{u_k})\,\d t \\
          &\quad = \limsup_{k
          \rightarrow \infty} \biggl( \eps \int_0^T
          (\ddot{y}_{\eps}^{u_k},y_{\eps}^{u_k}) \, \d t - \int_0^T
          (\dot{y}_{\eps}^{u_k},y_{\eps}^{u_k})\, \d t + \int_0^T
          (u_k,y_{\eps}^{u_k}) \, \d t	\biggr)\\
          & \quad = \limsup_{k
          \rightarrow \infty} \biggl( 	-\epsi
            (\dot{y}^{u_k}_\eps(0),y^0)  -  \epsi \int_0^T \|
             \dot{y}^{u_k}_\eps \|^2 \d t \\
          &\qquad -\frac{1}{2}\lVert
   y_{\eps}^{u_k}(T)\|^2 + \frac{1}{2}\lVert
   y^0\|^2 + \int_0^T
          (u_k,y_{\eps}^{u_k}) \, \d t	\biggr)
	\end{align*}
        where we also used the conditions  
        \eqref{eulerolagrangeepsilon1}-\eqref{eulerolagrangeepsilon3}.
        Owing to the above
        convergences we infer   
	\begin{align*}
	&\limsup_{k\rightarrow \infty} \int_0^T
   (\eta_{\eps}^{u_k},y_{\eps}^{u_k}) \, \d t \quad
          \\ &
               \quad \leq 	-\epsi
            (\dot{y} (0),y^0)  -   \epsi \int_0^T \|
            \dot{y}\|^2 \d t -\frac{1}{2}\lVert
   y (T)\|^2 + \frac{1}{2}\lVert
   y ^0\|^2 + \int_0^T
          (u,y) \, \d 
            t\\
          &\quad = \int_0^T
   (\eta,y) \, \d t.
	\end{align*}
	This implies that $\eta \in \partial \phi (y)$ via
        \eqref{eq:limsup2},  see \cite[Prop.~2.5, p.~27]{Brezis73}.  Hence, $y$ in the unique solution to the
        Euler-Lagrange problem
        \eqref{eulerolagrangeepsilon}-\eqref{eulerolagrangeepsilon3}
        and $y=y^u_\epsi$, $\eta=\eta^u_\epsi$ by Proposition
        \ref{ms13}.
      \end{proof}

\begin{lemma}[Coercivity of $P_{\eps\lambda}$]\label{coercivity}
We have
	\begin{equation}
	\label{coerc1}
		\eps^3 e^{-T/\eps} \lVert y - y_{\eps}^{u}
                \rVert_{H^1(0,T;H)}^2 \leq W_{\eps}(y,u) -
                M_\epsi^u\quad \forall (y,u)\in K(y^0)\times U.
	\end{equation}
	In particular, the sublevels of $\lambda \eps^{-3} e^{T/\eps}
        P_{\eps\lambda}(\cdot, u)$ are bounded in $H^1(0,T;H)$
        independently of $u\in U$.
\end{lemma}
\begin{proof}
	Let us start by rewriting $W_{\eps}(y,u) - M_\epsi^u$ as 
	\begin{align*}
          &W_{\eps}(y,u) - M_\epsi^u \\
          &\quad = \int_0^T e^{-t/\eps} \biggl(\frac{\eps}{2} \lVert
            \dot{y} \rVert^2 + \frac{\kappa}{2} \lVert y \rVert^2
            \biggr) \,\d t  + \left(\int_0^T e^{-t/\eps}\biggl( \psi(y) - (u,y)	\biggr) \,\d t - M_\epsi^u\right)\\
&\quad =: Q_{\eps}(y) + R_{\eps}(y,u),
	\end{align*}
	where the functional $R_{\eps}(\cdot,u)$ is convex and the
        quadratic functional $Q_{\eps}$ can be written in terms of $v = e^{-t/(2\eps)}y$ as 
	\begin{align*}
	Q_{\eps}(y) &= \int_0^T \biggl(	\frac{\eps}{2} \lVert \dot{v} \rVert^2 + \frac{1 + 4\eps \kappa}{8 \eps} \lVert v \rVert^2	\biggr) \,\d t + \frac{1}{4} \lVert v(T) \rVert^2 - \frac{1}{4} \lVert v(0) \rVert^2\\
	&=: V_{\eps}(v) + \frac{1}{4} \lVert v(T) \rVert^2 - \frac{1}{4} \lVert v(0) \rVert^2.
	\end{align*}
         We now use the fact that $V_\eps$ is quadratic. 
	For all $v_1,\, v_2 \in H^1(0,T;H)$ and all $r\in (0,1)$ one
         computes  
	\begin{align*}
          &V_{\eps}(rv_1 + (1-r)v_2) \nonumber\\
& =	\frac{\eps}{2} \int_0^T \| r\dot v_1 + (1-r)\dot v_2\|^2 \, \d
                                                   t + \frac{1  +   4  \eps \kappa}{8   \eps} \int_0^T \| r  v_1 + (1-r) v_2\|^2  \, \d  t\nonumber\\
& =   \frac{\eps }{2} \int_0^T\left(r^2\  \| \dot v_1\|^2 + (1-r)^2   \|  \dot v_2\|^2      +   2r      (1-r) (\dot  v_1, \dot v_2)  \right)\,\d t\nonumber\\
 & \quad+  \frac{1  +   4  \eps \kappa}{8   \eps}  \int_0^T\left(r^2\  \|  v_1\|^2 + (1-r)^2   \|   v_2\|^2      +   2r      (1-r) (  v_1,  v_2)  \right)\,\d t\nonumber\\
 & =   \frac{\eps }{2} \int_0^T\left(r  \| \dot v_1\|^2 + (1-r)
   \|  \dot v_2\|^2 - r (1-r) \|\dot  v_1- \dot v_2\|^2  \right)\,\d
   t\nonumber\\
          & \quad +    \frac{1   +   4  \eps \kappa}{8   \eps}
            \int_0^T\left(r  \|  v_1\|^2 + (1-r)   \|   v_2\|^2 - r
            (1-r) \|  v_1-  v_2\|^2  \right)\,\d t\nonumber\\   & =    r V_{\eps}(v_1) + (1 - r) V_{\eps}(v_2) - r(1 - r) \int_0^T \biggl(	\frac{\eps}{2} \lVert \dot{v}_1 - \dot{v}_2 \rVert^2 + \frac{1 +   4 \eps \kappa}{8 \eps} \lVert v_1 - v_2 \rVert^2	\biggr) \,\d t,
	\end{align*}
	which,  for  $\eps<\kappa$  small enough,  implies    the following inequality
	\begin{align}
	\label{convex1}
	V_{\eps}(rv_1 + (1-r)v_2) \leq r V_{\eps}(v_1) + (1 - r) V_{\eps}(v_2) - \eps \frac{r(1 - r)}{2} \lVert v_1 - v_2 \rVert_{H^1(0,T;H)}^2.
	\end{align}
        This in particular entails that 
	\begin{align}
	\label{convex2}
	W_{\eps}&(r y + (1 -r) y_{\eps}^{u} , u) -M_{\eps}^u \nonumber \\=& \, V_{\eps} (r v + (1 -r) v_{\eps}^{u}) + \frac{1}{4} \lVert r v(T) + (1 - r) v_{\eps}^{u}(T) \rVert^2 \nonumber \\ &-\frac{1}{4} \lVert r v(0) + (1 - r) v_{\eps}^{u}(0) \rVert^2 + R_{\eps}(r y + (1 - r) y_{\eps}^{u} ,u)\nonumber\\
	\leq& \, r V_{\eps}(v) + (1 - r)V_{\eps}(v_{\eps}^{u}) - \eps \frac{r(1 - r)}{2} \lVert v - v_{\eps}^{u} \rVert_{H^1(0,T;H)}^2\nonumber \\
	& + \frac{r}{4} \lVert v(T) \rVert^2 + \frac{1 - r}{4} \lVert v_{\eps}^{u}(T) \rVert^2 - \frac{1}{4} \lVert y^0 \rVert^2+ r R_{\eps}(y,u) + (1-r) R_{\eps}(y_{\eps}^{u},u) \nonumber\\
	=& \, r( W_{\eps}(y,u) - M_\epsi^u ) + (1-r)(W_{\eps}(y_{\eps}^u,u) - M_\epsi^u) - \eps \frac{r(1 - r)}{2} \lVert v - v_{\eps}^{u} \rVert_{H^1(0,T;H)}^2\nonumber\\
	=& \, r( W_{\eps}(y,u) - M_\epsi^u ) - \eps \frac{r(1 - r)}{2} \lVert v - v_{\eps}^{u} \rVert_{H^1(0,T;H)}^2,
	\end{align}
	where we used the convexity of the maps $v\mapsto
        R_{\eps}(v,u)$ and $v \mapsto \lVert v \rVert^2$, inequality
        \eqref{convex1}, the fact that $v(0) = v_{\eps}^{u}(0) = y^0$,
        and $W_{\eps}(y_{\eps}^{u},u) = M_\epsi^u$.  By observing that
        $$	W_{\eps} (r y + (1 -r) y_{\eps}^{u} , u) -M_{\eps}^u
        \geq\min_y W_{\eps}(y,u)-M_\epsi^u =   W_{\eps} ( y_{\eps}^{u} , u) -M_{\eps}^u=0,$$
        inequality \eqref{convex2} implies that
        $$r  ( W_{\eps}(y,u) - M_\epsi^u ) \geq  \eps \frac{r(1 -
          r)}{2} \lVert v - v_{\eps}^{u} \rVert_{H^1(0,T;H)}^2.$$
        Assume now that $r>0$, divide by $r$, and take $r\to
        0$ to get
        $$W_{\eps}(y,u) - M_\epsi^u \geq  \frac{\eps}{2}  \lVert v -
        v_{\eps}^{u} \rVert_{H^1(0,T;H)}^2.$$ 
	Eventually, putting $v= e^{-t/(2\eps)} y$ and $v_{\eps}^u =
        e^{-t/(2\eps)} y_{\eps}^u$  we obtain
	\begin{align*}
		W_{\eps}&(y,u) - M_\epsi^u \geq  \frac{\eps}{2}  \lVert v -
        v_{\eps}^{u} \rVert_{H^1(0,T;H)}^2\\
		\geq& \frac{\eps  }{2}   \int_0^T e^{-t/\eps} \biggl( \lVert y - y_{\eps}^u \rVert^2 + \lVert \dot{y} - \dot{y}_{\eps}^u \rVert^2 + \frac{\lVert y - y_{\eps}^u \rVert^2}{4 \eps^2} - \frac{1}{\eps} (y - y_{\eps}^u , \dot{y} - \dot{y}_{\eps}^u) 	\biggr) \, \d t \\
		\geq& \frac{\eps  }{2}   \int_0^T e^{-t/\eps} \biggl(	\lVert y - y_{\eps}^u \rVert^2 + \lVert \dot{y} - \dot{y}_{\eps}^u \rVert^2 + \frac{\lVert y - y_{\eps}^u \rVert^2}{4 \eps^2}   - \frac{1}{1 + 2\eps^2} \lVert \dot{y} - \dot{y}_{\eps}^u \rVert^2 \\&- \frac{1 + 2 \eps^2}{4 \eps^2} \lVert y - y_{\eps}^u \rVert^2 \biggr) \, \d t \\
		\geq& \,\eps^3 e^{-T/\eps} \lVert y - y_{\eps}^u \rVert_{H^1(0,T;H)}^2
	\end{align*}
        for $\varepsilon \in (0,1]$, 
	which proves \eqref{coerc1}. In particular, we have
	\begin{align*}
		\lVert y \rVert_{H^1(0,T;H)}^2 \leq& 2 \lVert y - y_{\eps}^u \rVert_{H^1(0,T;H)}^2 + 2 \lVert y_{\eps}^u \rVert_{H^1(0,T;H)}^2\\
		\leq& c + 2 \lambda \eps^{-3} e^{T/\eps} P_{\eps\lambda} (y,u),
	\end{align*}
	which concludes the proof.
      \end{proof}

\begin{lemma}
	\label{argument1}
	Let $u_{\eps} \in U$ with $u_{\eps}\rightarrow u$ in
        $L^2(0,T;H)$. Then, up to subsequences $y_{\eps}^{u_\eps} \rightharpoonup y$ in $H^1(0,T;H)$ and $y_{\eps}^{u_\eps}  \rightarrow y$ in $H^{\sigma}(0,T;H)$ for any $\sigma \in (0,1)$ with $y=S(u)$.
\end{lemma}
\begin{proof}
  Letting $\eta_\epsi = u_\epsi +\epsi\ddot{y}^{u_\epsi}_\epsi - \dot{y}^{u_\epsi}_\epsi $,
the uniform bound \eqref{maxreg} and up to not relabeled subsequences
we have that 
	\begin{align}
	\eps\ddot{y}_{\eps}^{u_\epsi} \rightharpoonup 0 \quad \text{in}\; L^2(0,T;H),\label{conveps1}\\
	 y_{\eps}^{u_\eps} \rightharpoonup y \quad \text{in}\; H^1(0,T;H),\label{conveps2}\\
	\eta_{\eps}^{u_\epsi} \rightharpoonup \tilde{\eta} \quad \text{in} \; L^2(0,T;H),\label{conveps3}
	\end{align}
	for some functions $y$ and $\tilde{\eta}$ with $
	\dot{y} + \tilde{\eta} = u $.
	
We shall now check that $y_{\eps}^{u_\epsi}$ is indeed a Cauchy
sequence in $C^0([0,T];H)$. Let $y_{\eps}^{u_\epsi}$ and
$y_{\mu}^{u_\mu}$ be solutions of the problem
\eqref{eulerolagrangeepsilon}-\eqref{eulerolagrangeepsilon3} at level
$\eps$ and $\mu$, respectively. Consider the difference of the
Euler-Lagrange equation \eqref{eulerolagrangeepsilon} at level $\eps$
and the one at level $\mu$ and take its scalar product with the
function $w := y_{\eps}^{u_\epsi} - y_{\mu}^{u_\mu}$. By letting
$\eta_\mu=u_\mu +\mu\ddot{y}^{u_\mu}_\mu - \dot{y}^{u_\mu}_\mu $ and integrating in time over $(0,t)$ we get
\begin{align*}
	&- \int_0^t (\eps \ddot{y}_{\eps}^{u_\epsi} - \mu
  \ddot{y}_{\mu}^{u_\mu}, w) \, \d t + \int_0^t \frac{\d}{\d
  t}\frac{1}{2} \lVert w \rVert^2 \, \d t + \int_0^t (\eta_{\eps} -
  \eta_{\mu} , y_{\eps}^{u_\epsi} - y_{\mu}^{u_\mu}) \, \d t \\
  &\quad = \int_0^t (u_{\eps}-u_{\mu} , w)\, \d t.
\end{align*}
The $\kappa$-convexity of $\phi$, the fact that $w(0) = 0$, and an integration by parts give
\begin{align*}
	&\eps \int_0^t \lVert \dot{w} \rVert^2 \, \d t +
   \kappa \int_0^t \lVert w \rVert^2 \, \d t +
   \frac{1}{2}\lVert w(t) \rVert^2 \\
        &\quad \leq (\eps \dot{y}_{\eps}^{u_\epsi}(t) - \mu
   \dot{y}_{\mu}^{u_\mu} (t) , w(t)) - (\eps -\mu) \int_0^t (\dot{y}_{\mu}^{u_\mu} ,\dot{w}) \, \d t + \int_0^t (u_{\eps} - u_{\mu} , w) \, \d t.
\end{align*}
Young's inequality allows then to deduce
\begin{align}
	&\eps \int_0^t \lVert \dot{w}\rVert^2 \, \d t +
   \kappa\int_0^t \lVert w \rVert^2 \, \d t +
   \frac{1}{ 4}\lVert w (t) \rVert^2\nonumber \\
  &\quad \leq  2  \eps^2 \lVert \dot{y}_{\eps}^{u_\eps}
    \rVert_{C^0([0,T];H)}^2 + 2 \mu^2 \lVert \dot{y}_{\mu}^{u_\mu} \rVert_{C^0([0,T];H)}^2 \nonumber \\[1mm] &\quad + (\eps + \mu) \lVert
                                           \dot{y}_{\mu}^{u_\mu}
                                           \rVert_{L^2(0,T;H)} \lVert
                                           \dot{w} \rVert_{L^2(0,T;H)}
                                           + \lVert u_{\eps} - u_{\mu} \rVert_{L^2(0,T;H)} \lVert w \rVert_{L^2(0,T;H)}.\label{ineq}
\end{align}
The Gagliardo-Nirenberg inequality \cite[Comments $(iii)$,
p. 233]{brezis}, \cite[Theorem 1, p. 734]{niremberg66}, and bound \eqref{maxreg} give
\begin{align*}
	&\lVert \dot{y}_{\eps}^{u_\eps} \rVert_{C^0([0,T];H)} \leq c \lVert
        \ddot{y}_{\eps}^{u_\eps}\rVert_{L^2(0,T;H)}^{1/2} \lVert
        \dot{y}_{\eps}^{u_\eps} \rVert_{L^2(0,T;H)}^{1/2} + c \lVert
        \dot{y}_{\eps}^{u_\eps} \rVert_{L^2(0,T;H)} \leq c
   \eps^{-1/2},\\
  &
   \lVert \dot{y}_{\mu}^{u_\mu} \rVert_{C^0([0,T];H)} \leq c \lVert
        \ddot{y}_{\mu}^{u_\mu}\rVert_{L^2(0,T;H)}^{1/2} \lVert
        \dot{y}_{\mu}^{u_\mu} \rVert_{L^2(0,T;H)}^{1/2} + c \lVert
        \dot{y}_{\mu}^{u_\mu} \rVert_{L^2(0,T;H)} \leq c \mu^{-1/2}.
\end{align*}
We can hence use \eqref{ineq}  and the fact that $\lVert
\dot{y}_{\mu}^{u_\mu}   \rVert_{L^2(0,T;H)}\leq c$  and $\lVert \dot{w}
    \rVert_{L^2(0,T;H)}\leq c$  to obtain
\begin{align*}
\eps \int_0^t \lVert \dot{w}\rVert^2 \, \d t + \frac{1}{ 4}\lVert w (t)
  \rVert^2 \leq c (\eps + \mu) + c \lVert u_{\eps} - u_{\mu}
  \rVert_{L^2(0,T;H)}  - \kappa \int_0^t \| w\|^2\, \d t.
\end{align*} 
 This entails that 
\begin{equation}
	\label{gronw}
	\lVert y_{\eps} - y_{\mu} \rVert_{C^0([0,T];H)} \leq c \biggl( \eps + \mu + \lVert u_{\eps} - u_{\mu} \rVert_{L^2(0,T;H)} \biggr)^{1/2}
      \end{equation}
       where we also applied the Gronwall Lemma in case $\kappa<0$. The  strong convergence
\begin{equation}
\label{convStrong}
y_{\eps} \rightarrow y \;\; \text{in} \;C^0([0,T];H),
\end{equation}
follows.

To identify the limit function $\tilde{\eta}$, we  use convergences
\eqref{conveps3} and \eqref{convStrong} to get
\begin{align*}
	\limsup_{\eps \rightarrow 0} \int_0^T (\eta_{\eps},y_{\eps}) \, \d t = \int_0^T (\tilde{\eta},y)\, \d t.
\end{align*}
This implies that $\tilde \eta \in \partial \phi(u)$ a.e. by
\eqref{eq:limsup2}. The initial condition \eqref{gradientflowf}
follows from \eqref{convStrong}, so that $y =S(u)$. As this limit is
unique, the whole sequence $(y_\eps^{u_\epsi})_\eps$ converges and no
extraction of a subsequence is actually necessary. 

  We now make use of the interpolation space $
(C^0([0,T];H),H^1(0,T;H))_{\sigma,1} $  for $\sigma \in (0,1)$,  whose elements are 
those functions $w \in C^0([0,T];H)$ such that
$$\| w \|_{ (C^0([0,T];H),H^1(0,T;H))_{\sigma,1} } := \int_0^\infty
r^{-\sigma-1} K(r,w)\, \d r <\infty,$$
where   $K:(0,\infty) \times C^0([0,T];H) \to [0,\infty)$ is
defined as
\begin{align*}
  K(r,w):= &\inf\Big\{ \| w_0 \|_{C^0([0,T];H)} + r \|
    w_1\|_{H^1(0,T;H)}  \ : \ \\
    &\qquad w_0\in C^0([0,T];H), \ w_1\in H^1(0,T;H),
    \ w=w_0 +w_1 \Big\}.
    \end{align*}
see
the classical reference 
\cite{BL}. We will also use that there exists $c>0$ such
that 
$$ \| w \|_{ (C^0([0,T];H),H^1(0,T;H))_{\sigma,1} } \leq  c\lVert w
                 \rVert_{C^0([0,T];H)}^{1-\sigma}  \lVert
              w
                 \rVert_{H^1(0,T;H)}^{\sigma}  \quad \forall w \in
                 H^1(0,T;H),$$
                 see \cite[Lemma~2.1.i]{Moiola}.
As $ (C^0([0,T];H),H^1(0,T;H))_{\sigma,1} \subset H^{\sigma}(0,T;H)$
\cite[Theorem~6.2.4, p.~142]{BL},  the uniform bound of $y_{\eps}^{u_\eps}$ in $H^1(0,T;H)$ and estimate \eqref{gronw} imply that $y_{\eps}^{u_\eps} \rightarrow y$ in $H^{\sigma}(0,T;H)$. Indeed, one has
\begin{align*}
  &  \lVert y_{\eps}^{u_\eps} -
    y_{\mu}^{u_\mu}\rVert_{H^{\sigma}(0,T;H)}    \leq c  \lVert y_{\eps}^{u_\eps} - y_{\mu}^{u_\mu}
    \rVert_{(C^0([0,T];H),H^1(0,T;H))_{\sigma,1}} \\[2mm]
  &\quad \leq  c
                 \lVert y_{\eps}^{u_\eps} - y_{\mu}^{u_\mu}
                 \rVert_{C^0([0,T];H)}^{1-\sigma}  \lVert
                 y_{\eps}^{u_\eps} - y_{\mu}^{u_\mu}
                 \rVert_{H^1(0,T;H)}^{\sigma} \\ &\quad\leq c \biggl( \eps + \mu + \lVert u_{\eps} - u_{\mu} \rVert_{L^2(0,T;H)} \biggr)^{(1 - \sigma)/2} \rightarrow 0.\qedhere
\end{align*}
\end{proof}
\begin{lemma}
	\label{argument2}
	Let $(y_{\eps\lambda},u_{\eps\lambda}) \in H^1(0,T;H)\times U$
        be such that $P_{\eps\lambda}(y_{\eps\lambda},u_{\eps\lambda})
        \leq c$, with $c$ independent of $\lambda$. As
        $\lambda\rightarrow 0$, up to not relabeled subsequences we
        have that $(y_{\eps\lambda},u_{\eps\lambda})\overset{\tau}{\rightarrow} (y_{\eps}^{u_{\eps}},u_{\eps})$.
\end{lemma}
\begin{proof}
From the compact injection $U \subset \subset L^2(0,T;H)$ we have
$u_{\eps\lambda} \rightarrow u_{\eps}$ in $L^2(0,T;H)$ along some not
relabeled subsequence. Lemma \ref{coercivity} implies that
$(y_{\eps\lambda})_\lambda$ is  bounded in $H^1(0,T;H)$. Hence, up to not relabeled subsequences, $y_{\eps\lambda}\rightharpoonup y_{\eps}$ in $H^1(0,T;H)$.
	
As  $P_{\eps\lambda}(y_{\eps\lambda},u_{\eps\lambda})<c$  and $J$ is
nonnegative we have $$ 0 \leq
W_{\eps}(y_{\eps\lambda},u_{\eps\lambda}) - m_{\eps}(u_{\eps\lambda})
\leq c \lambda. $$ Recall from Section \ref{section6} that
$(y,u)\mapsto W_{\eps}(y,u) - M_{\eps}^u$ is lower semicontinuous with
respect to the topology $\tau$. Hence, by passing to the $\liminf$ as
$\lambda \rightarrow 0$, we get $$ W_{\eps}(y_{\eps},u_{\eps}) =
M_{\eps}^{u_{\eps}}, $$ proving indeed that $y_{\eps} =
y_{\eps}^{u_{\eps}}$, which concludes the argument.
\end{proof}

\begin{lemma}
	\label{argument2'}
	Let $\lambda = \lambda_\eps$, with $\limsup_{\eps \rightarrow
          0} \lambda_\eps\eps^{-3} e^{T/\eps} = 0$ and let $(y_{\eps\lambda},u_{\eps\lambda}) \in H^1(0,T;H)$  $\times U$  be such that $P_{\eps\lambda}(y_{\eps\lambda},u_{\eps\lambda}) \leq c$, with $c$ independent of $\eps$. Then, as $\eps \rightarrow 0$, up to not relabeled subsequences $(y_{\eps\lambda},u_{\eps\lambda}) \overset{\tau}{\rightarrow} (y,u)$ with $y=S(u)$. 
\end{lemma}
\begin{proof}
	Let $(y_{\eps\lambda},u_{\eps\lambda}) \in H^1(0,T;H)\times U$ be such that $P_{\eps\lambda}(y_{\eps\lambda},u_{\eps\lambda}) \leq c$.
	As the injection $U \subset \subset L^2(0,T;H)$ is compact, we
        get $u_{\eps\lambda} \rightarrow u$ in $L^2(0,T;H)$ along some
        not relabeled subsequence. Using Lemma \ref{coercivity}, we have
	\begin{equation}
	\label{coer}
	\limsup_{\eps \rightarrow 0} \lVert y_{\eps\lambda} - y_{\eps}^{u_{\eps\lambda}} \rVert_{H^1(0,T;H)}^2 \leq \limsup_{\eps \rightarrow 0} \lambda_\eps \eps^{-3} e^{T/\eps} P_{\eps\lambda} (y_{\eps\lambda},u_{\eps\lambda}) =0.
      \end{equation}
      On the other hand, Lemma \ref{argument1} implies that
      $y^{u_{\eps\lambda}}_\eps \rightharpoonup y$ in
      $H^1(0,T;H)$  and  $y^{u_{\eps\lambda}}_\eps \to y$ in
      $H^\sigma(0,T;H)$ for all $\sigma\in (0,1)$, with $y=S(u)$. This
      implies that $y_{\eps\lambda} \rightharpoonup y$ in $H^1(0,T;H)$.
	
	The thesis follows then from \eqref{coer} and Lemma
        \eqref{argument1} by simply observing that $$\lVert
        y_{\eps\lambda} - y \rVert_{H^{\sigma}(0,T;H)} \leq c \lVert y_{\eps\lambda} - y^{u_{\eps\lambda}}_{\eps} \rVert_{H^1(0,T;H)} + \lVert y^{u_{\eps\lambda}}_{\eps} - y \rVert_{H^{\sigma}(0,T;H)} \rightarrow 0 .\qedhere$$
\end{proof}


\section{Existence: proof of Theorem \ref{theorem}.{\rm i}}
\label{well}
In this section, we existence for Problems \eqref{initialproblem},
\eqref{epsilonproblem}, and \eqref{epsilonlambdaproblem}, namely
Theorem \ref{theorem}.i.

\subsection{Well-posedness for Problem \ref{initialproblem}}
\label{refP}
Taking any $\tilde{u} \in U$ letting $\tilde{y} = S(\tilde{u})$ one
has that $0\leq P(\tilde{y},\tilde{u}) = J(\tilde{y},\tilde{u}) <
\infty$. In particular $\inf_{H^1(0,T;H)\times U}  P \in [0,\infty)$.

Let $(y_k,u_k)_k \subset H^1(0,T;H)\times U$ be an infimizing sequence
for $P$, that is $P(y_k,u_k) \rightarrow \inf_{H^1(0,T;H)\times U}
P$. The strong convergence
\begin{equation}\label{conv0} u_k \rightarrow u \qquad \text{in} \;
  L^2(0,T;H) \end{equation} 
follows from the compact injection $U \subset \subset L^2(0,T;H)$, up to some not relabeled subsequence.
For all $k>0$, $y_k\in S(u_k) $ solves the gradient flow
\begin{align}
\dot{y}_k(t) + \eta_k(t) &= u_k(t) \quad \text{in} \; H, \quad \text{a.e.}\;t\in (0,T),\label{gradientflowk}\\
\eta_k(t) &\in \partial \phi(y_k(t)) \quad \text{in} \; H, \quad\text{a.e.}\;t\in (0,T),\\
y_k(0)&=y^0.\label{gradientflowk3}
\end{align}
As each term in equation \eqref{gradientflowk} is in $L^2(0,T;H)$, we
use the chain rule \eqref{eq:chain} in order to compute
\begin{align*}
&\int_0^T \lVert \dot{y}_k \rVert^2 \,\d t+ \int_0^T \lVert \eta_k
  \rVert^2 \,\d t = \int_0^T \lVert \dot{y}_k + \eta_k \rVert^2 \,\d t
  - 2 \int_0^T (\dot{y}_k,\eta_k)\,\d t \\ &\quad = \int_0^T \lVert
                                             u_k \rVert^2 \,\d t- 2
                                             \int_0^T \frac{\d}{\d t}
                                             \phi(y_k) \,\d t = \int_0^T \lVert
                                             u_k \rVert^2 \,\d t-
                                              2\phi(y_k(T))  {}+ 2
                                             \phi(y^0)\\
  &\quad \leq \int_0^T \lVert
                                             u_k \rVert^2 \,\d t-
                                             2  ((\partial
    \phi(y^0))^\circ,  y_k (T) {} - y^0)\\
  &\quad \leq \int_0^T \lVert
                                             u_k \rVert^2 \,\d t
    +\frac12\int_0^T \| \dot{y}_k\|^2  \d t  {} + 2T \| (\partial \phi(y^0))^\circ\|^2
\end{align*}
where we used the equation \eqref{gradientflowk} as well as the fact
that $y^0 \in D(\partial \phi)$ . We have obtained $\lVert \dot{y}_k \rVert_{L^2(0,T;H)} + \lVert \eta_k \rVert_{L^2(0,T;H)} \leq c$, which yields, up to not relabeled subsequences, that
\begin{align}
&y_k \rightharpoonup y \qquad \text{in}\; H^1(0,T;H),\label{conv1}\\
&\eta_k \rightharpoonup \eta \qquad \text{in} \; L^2(0,T;H),
\end{align}
for some limit functions $y,\eta$ with  $ \dot{y} + \eta = u $  a.e.   As $y_k(0) \rightharpoonup y(0)$, we have that $y(0)=y^0$. Moreover,
\begin{align*}
	&\limsup_{k \rightarrow \infty} \int_0^T (\eta_k, y_k) \,\d t = \limsup_{k \rightarrow \infty} \int_0^T (u_k - \dot{y}_k, y_k) \,\d t
	\\ &\quad =\int_0^T (u,y) \,\d t - \liminf_{k\rightarrow \infty} \frac{1}{2}\lVert y_k(T) \rVert^2 + \frac{1}{2}\lVert y^0 \rVert^2\\
	&\quad \leq  \int_0^T (u - \dot{y} , y) \,\d t =  \int_0^T (\eta, y) \,\d t.
\end{align*}
Again \eqref{eq:limsup2} entails that $\eta \in \partial \phi(y)$ in
$L^2(0,T;H)$. The limit function $y$ hence solves then the gradient flow problem, namely, $y \in S(u)$. Eventually, convergences \eqref{conv0} and \eqref{conv1} together with Assumption (A3) imply
\begin{equation*}
\inf_{H^1(0,T;H)\times U} P \leq  J(y,u) \leq \liminf_{k\rightarrow \infty} J(y_k, u_k) = \inf_{H^1(0,T;H)\times U} P ,
\end{equation*}
so that $(y,u)$ actually solves \eqref{initialproblem}.

\subsection{Well-posedness for Problem \ref{epsilonproblem}}
\label{Pepsilonadmitssolution}

Choosing an arbitrary $\tilde u \in U $ and letting $\tilde y \in
S(\tilde u)$ one has that $0\leq P_\epsi(\tilde y \, \tilde u) =
J(\tilde y , \tilde u)<\infty$. In particular, $\inf_{H^1(0,T;H)\times
   U } P_{\eps} \in [0,\infty)$.

	Let $(y_k,u_k)_k \subset H^1(0,T;H) \times U$ be an infimizing
        sequence for $P_{\eps}$, namely, $	P_{\eps}(y_k,u_k)
        \rightarrow \inf_{H^1(0,T;H)\times U} P_{\eps}$. We have,
        $P_{\eps}(y_k,u_k) = J(y_k,u_k)$ and $y_k \in \arg\min
        W_{\eps}(\cdot,u_k)$, namely, $y_k = y_{\eps}^{u_k}$. The
        strong convergence  $u_k \to u$ in $L^2(0,T;H)$ follows by the compact injection $U \subset \subset L^2(0,T;H)$, up to a not relabeled subsequence. For each $k>0$, there exists $\eta_k \in L^2(0,T;H)$ such that
	\begin{align*}
	-\eps \ddot{y}_k(t) + \dot{y}_k(t) + \eta_k(t) &= u_k(t) \quad
                                                         \text{for
                                                         a.e.}\ t\in(0,T),\\
	\eta_k(t) &\in \partial\phi(y_k(t)) \quad  \text{for a.e.}\ t\in(0,T),\\
	y_k(0) &= y^0,\\
	\eps\dot{y}_k(T) &= 0.
	\end{align*}
Lemma \ref{lemma2} ensures that, up to not relabeled subsequences,
$y_k\rightharpoonup y$ in $H^2(0,T;H)$, $\eta_k\rightharpoonup \eta$
in $L^2(0,T;H)$ where $\eta \in \partial \phi(y)$ a.e., and $y$ solves
the Euler-Lagrange problem \eqref{eulerolagrangeepsilon}-\eqref{eulerolagrangeepsilon3} corresponding to $u$. In particular, $y =y^u_\eps$. Together with Assumption (A3), these convergences imply
	\begin{equation*}
\inf_{H^1(0,T;H)\times U} P_{\eps} \leq J(y,u)\leq \liminf_{k\rightarrow \infty} J(y_k,u_k) = \inf_{H^1(0,T;H)\times U} P_{\eps}
	\end{equation*}
	which proves that $(y,u)$ actually minimizes $P_{\eps}$.

\subsection{Well-posedness for Problem \ref{epsilonlambdaproblem}}
\label{section6}

Given any $\tilde u \in U$ we have that $0 \leq
P_{\epsi\lambda}(y^{\tilde u}_\epsi,\tilde u) $  $= J (y^{\tilde
  u}_\epsi,\tilde u)<\infty$. This proves that $\inf_{H^1(0,T;H)\times
  U} P_{\eps\lambda}\in [0,\infty)$.

Let $(y_k,u_k)_k \subset H^1(0,T;H)\times U$ be an infimizing sequence
for $P_{\eps\lambda}$, namely, such that $P_{\eps\lambda}(y_k,u_k) \rightarrow
\inf_{H^1(0,T;H)\times U} P_{\eps\lambda}$.    The strong convergence
$u_k \to u$ in  $L^2(0,T;H)$ follows from the compact injection  $U
\subset \subset L^2(0,T;H)$, up to a not relabeled subsequence.
Using Lemma \ref{coercivity}, we have that $\lVert y_k \rVert_{H^1(0,T;H)} \leq c$. Then, up to not relabeled subsequences, the following convergence hold
\begin{align}
y_k \rightharpoonup y \qquad \text{in}\; H^1(0,T;H).\label{convergence2}
\end{align}
Since $H^1(0,T;H) \subset C^0([0,T];H)$, we have that $y_k(t)
\weakto y(t)$ for all $t\in[0,T]$. In particular, the initial condition $y(0)=y^0$ is satisfied.

Lemma \ref{lemma2} implies that
\begin{equation}
y_{\eps}^{u_k} \rightharpoonup y_{\eps}^u \qquad \text{in}\; H^2(0,T;H) \label{convergence3},
\end{equation}
and $y_{\eps}^u$ satisfies the Euler-Lagrange problem \eqref{eulerolagrangeepsilon}-\eqref{eulerolagrangeepsilon3} corresponding to $u$. 
By Assumption (A3), we have that $J(y,u) \leq \liminf_n J(y_n,u_n)$.

We hence reduce ourselves to check that $(y,u)\mapsto W_{\eps}(y,u) -
M_{\eps}^u$ is lower semicontinuous with respect to the topology $\tau$. By using Lemma \ref{lemma1}, we have
\begin{align*}
&W_{\eps} (y_k,u_k) - M_{\eps}^{u_k}\\ &\quad = \int_0^T e^{-t/\eps}
                                         \biggl( \frac{\eps}{2} \lVert
                                         \dot{y}_k \rVert^2 +
                                         \phi(y_k) -(y_k,u_k)
                                         \biggr) \,\d t  +
                                         \frac{\eps^2}{2} \lVert
                                         \dot{y}_{\eps}^{u_k}(0)\rVert^2
  \\ &\qquad+ \eps e^{-T/\eps} \phi(y_{\eps}^{u_k}(T)) - \eps \phi(y^0) - \eps \int_0^T e^{-t/\eps} (u_k , \dot{y}_{\eps}^{u_k}) \,\d t + \int_0^T e^{-t/\eps} (u_k,y_{\eps}^{u_k}) \,\d t.
\end{align*}
Taking the {$\liminf$} as $k \rightarrow \infty$ and using
convergences \eqref{conv0} and \eqref{convergence3}, the lower semicontinuity of $\lVert \cdot \rVert$ and of $\phi$, and the fact that we have $\dot{y}_{\eps}^{u_k}(t)\rightharpoonup \dot{y}_{\eps}^{u}(t)$ in $H$ for every $t \in [0,T]$ from convergence \eqref{convergence3}, we obtain
\begin{align*}
&\liminf_{k \rightarrow \infty} \biggl( W_{\eps} (y_k,u_k)-
                 M_{\eps}^{u_k} \biggr)\\
  &\quad \geq\int_0^T e^{-t/\eps} \biggl( \frac{\eps}{2}\lVert \dot{y}
    \rVert^2 + \phi(y) 	-(y,u)	\biggr) \,\d t + \frac{\eps^2}{2}
    \lVert \dot{y}_{\eps}^{u}(0)\rVert^2\\ &
                                             \qquad+ \eps e^{-T/\eps} \phi(y_{\eps}^{u}(T)) - \eps \phi(y^0) - \eps \int_0^T e^{-t/\eps} (u , \dot{y}_{\eps}^{u}) \,\d t+ \int_0^T e^{-t/\eps} (u,y_{\eps}^{u}) \,\d t\\
&\quad = W_{\eps} (y,u) - M_\epsi^u.
\end{align*}
We have checked that $P_{\epsi\lambda}(y,u)\leq \liminf_k
P_{\epsi\lambda}(y_k,u_k) = \inf_{H^1(0,T;H)\times U} P_{\eps\lambda}$, proving that 
$(y,u)$ minimizes $P_{\eps\lambda}$.

\section{Convergence: proof of Theorem \ref{theorem}.{\rm ii-iv}}
\label{proof}

\subsection{Proof of Theorem \ref{theorem}.ii.}
\label{52}

Let us start by checking the $\Gamma$-convergence $P_\eps
\stackrel{\Gamma}{\to} P$. We focus first on the $\Gamma$-$\liminf$  
inequality \eqref{liminf1}. Assume $(y_{\eps},u_{\eps})
\overset{\tau}{\rightarrow} (y,u)$. Without loss of generality,
$\sup_\epsi P_{\eps}(y_{\eps},u_{\eps}) <\infty$. Then, $P_{\eps} (y_{\eps},u_{\eps}) = J(y_{\eps},u_{\eps})$, $y_{\eps} = y^{u_{\eps}}_\epsi$, and we have
	\begin{equation*}
	\liminf_{\eps \rightarrow 0} P_{\eps}(y_{\eps},u_{\eps}) = \liminf_{\eps \rightarrow 0} J(y_{\eps},u_{\eps}) \geq J(y,u),
	\end{equation*}
	where we used the lower semicontinuity of $J$. The identification $J(y,u) = P(y,u)$, and hence the inequality \eqref{liminf1}, follows then from Lemma \ref{argument1}.
	
	To prove the recovery-sequence condition \eqref{limsup1}, we
        can assume, without loss of generality, that
        $P(\hat{y},\hat{u}) < \infty$. Then, we have $\hat{y} =
        S(\hat{u})$. As a recovery sequence we choose
        $\hat{u}_{\eps}=\hat{u}$ and $\hat{y}_{\eps} =
        y_{\eps}^{\hat{u}} $. Then, we have the identification
        $P_{\eps}(\hat{y}_{\eps},\hat{u}_{\eps}) =
        J(\hat{y}_{\eps},\hat{u})$. Moreover,   Lemma \ref{argument1}
        ensures convergence $\hat{y}_{\eps}\rightarrow \hat{y}$ in
        $H^{\sigma}(0,T;H)$ for all $\sigma\in(0,1)$. Assumption (A3)
        entails that 
	\begin{equation*}
	\limsup_{\eps \rightarrow 0}
        P_{\eps}(\hat{y}_{\eps},\hat{u}_{\eps}) = \limsup_{\eps
          \rightarrow 0} J(\hat{y}_{\eps},\hat{u}) \leq J(\hat{y},\hat{u}),
	\end{equation*}
	where we used the upper semicontinuity of $J(\cdot,\hat{u})$
        in the strong topology of $H^{\sigma}(0,T;H)$. The
        identification $J(\hat{y},\hat{u})=P(\hat{y},\hat{u})$ follows
        again from Lemma \ref{argument1}, which in particular ensures
        that $\hat{y}=S(\hat{u})$.
        
As the functionals $P_{\eps}$ are equicoercive in $H^1(0,T;H) \times
L^2(0,T;H)$ from Proposition \ref{ms13}, Theorem \ref{theorem}.ii
follows from the Fundamental Theorem of $\Gamma$-convergence \cite[Thm.~7.4, p.~69]{dalMaso}.

\subsection{Proof of Theorem \ref{theorem}.iii.}
\label{51}

In order to prove the $\Gamma$-convergence $P_{\eps\lambda}
\overset{\Gamma}{\rightarrow} P_{\eps}$, let us first check the
 $\Gamma$-$\liminf$ inequality \eqref{liminf1}. Without loss of
generality, let $(y_{\eps\lambda},u_{\eps\lambda})
\overset{\tau}{\rightarrow} (y_{\eps},u_{\eps})$ as $\lambda
\rightarrow 0$ be such that $\sup_\lambda P_{\eps\lambda}(y_{\eps\lambda},u_{\eps\lambda}) <\infty$. Then, we get $$ \liminf_{\lambda \rightarrow 0} P_{\eps\lambda}(y_{\eps\lambda},u_{\eps\lambda})   \geq \liminf_{\lambda \rightarrow 0} J(y_{\eps\lambda},u_{\eps\lambda}) \geq J(y_{\eps},u_{\eps}), $$ where the last inequality follows from the lower semicontinuity of $J$. Eventually, the identification $J(y_{\eps},u_{\eps}) = P_{\eps}(y_{\eps},u_{\eps})$ directly follows from Lemma \ref{argument2}.

A recovery sequence is given by 
$(\hat{y}_{\eps\lambda},\hat{u}_{\eps\lambda})=(y_{\eps}^{\hat{u}_{\eps}},\hat{u}_{\eps})$. In
  fact,
we readily obtain $$ \limsup_{\lambda\rightarrow 0}
P_{\eps\lambda}(y_{\eps}^{\hat{u}_{\eps}},\hat{u}_{\eps\lambda}) =
\limsup_{\lambda \rightarrow 0} J(y_{\eps}^{\hat{u}_{\eps}},\hat{u}_{\eps}) =
P_{\eps}(\hat{y}_{\eps},\hat{u}_{\eps}). $$

The equicoerciveness of the functionals $P_{\eps\lambda}$ in
$H^1(0,T;H) \times L^2(0,T;H)$ from Lemma \ref{coercivity} and an
application of the Fundamental Theorem of $\Gamma$-convergence
\cite[Thm.~7.4, p.~69]{dalMaso} conclude the proof.

\subsection{Proof of Theorem \ref{theorem}.iv.}

Recall that $\lambda = \lambda_\eps$ is such that $\limsup_{\eps
  \rightarrow 0} \lambda_\eps \eps^{-3} e^{T/\eps} =0$. In order to
check that $P_{\eps\lambda_\epsi} \overset{\Gamma}{\rightarrow} P$ we start
from the $\Gamma$-$\liminf$ inequality \eqref{liminf1}. Without loss of
generality, we can consider $\sup_\epsi P_{\eps\lambda_\eps}(y_{\eps\lambda_\epsi},u_{\eps\lambda_\epsi}) <\infty $. Hence, from the definition of $P_{\eps\lambda_\epsi}$, we deduce
	\begin{equation*}
		\liminf_{\eps \rightarrow 0} P_{\eps\lambda_\epsi} (y_{\eps\lambda_\epsi},u_{\eps\lambda_\epsi})  \geq \liminf_{\eps \rightarrow 0} J(y_{\eps\lambda_\epsi},u_{\eps\lambda_\epsi}) \geq J(y,u),
	\end{equation*}
	where we used the fact that
        $W_{\eps}(y_{\eps\lambda_\epsi},u_{\eps\lambda_\epsi}) -
        M_{\eps}^{u_{\eps\lambda_\epsi}} \geq 0$ and the lower semicontinuity of $J$. The identification $J(y,u)=P(y,u)$ follows then from Lemma \ref{argument2'}.
	
As regards the recovery-sequence condition \eqref{limsup1}, we first
note that, without loss of generality, one can assume
$P(\hat{y},\hat{u}) <\infty$, which implies $\hat{y} = S(\hat{u})$. We
choose the recovery sequence $(y_{\eps}^{\hat{u}},\hat{u})$, so that
$P_{\eps\lambda_\epsi}((y_{\eps}^{\hat{u}},\hat{u}) =
J(y_{\eps}^{\hat{u}},\hat{u})$. Lemma \ref{argument1} implies the
convergence $y_{\eps}^{\hat{u}} \rightarrow \hat{y}$ in
$H^{\sigma}(0,T;H)$ for all $\sigma\in(0,1)$. By exploiting the upper
semicontinuity of $J(\cdot,\hat{u})$ in the strong $H^{\sigma}(0,T;H)$
topology for some $\sigma \in (0,1)$ from assumption (A3), this
entails that 
	\begin{equation*}
	\limsup_{\eps  \rightarrow 0}
        P_{\eps\lambda_\epsi}(y_{\eps}^{\hat{u}},\hat{u}) =
        \limsup_{\eps \rightarrow 0} J(y_{\eps}^{\hat{u}},\hat{u}) \leq J(\hat{y},\hat{u}).
    \end{equation*}
    As we have that $\hat{y}=S(\hat{u})$, the equality
  $J(\hat{y},\hat{u})=P(\hat{y},\hat{u})$ follows.

Having proved the $\Gamma$ convergence, the statement follows from the 
equicoerciveness of the functionals $P_{\eps\lambda_\epsi}$ in
$H^1(0,T;H) \times L^2(0,T;H)$ from Lemma \ref{coercivity} by applying
again \cite[Thm.~7.4, p.~69]{dalMaso}.

\section*{Acknowledgements}
 T.F. acknowledges the support from the JSPS KAKENHI Grant-in-Aid for Scientific Research(C), Japan, Grant Number 21K03309.  
U.S. is supported by the Austrian Science Fund (FWF) through
projects 10.55776/F65, 10.55776/I4354, 10.55776/I5149, and 10.55776/P32788.

\end{document}